\documentclass[a4paper,12pt]{article}
\usepackage{amsfonts,amssymb,graphicx}

\newtheorem{theorem}{Theorem}

\begin{document}

\author{
Octavian G. Mustafa\\
\small{University of Vienna,}\\
\small{Faculty of Mathematics,}\\
\small{Nordbergstrasse 15, A-1090 Vienna, Austria}\\
\small{e-mail address: octawian@yahoo.com}
}

\title{On the uniqueness of flow in a recent tsunami model}
\date{}
\maketitle

\noindent{\bf Abstract} We give an elementary proof of uniqueness for the integral curve starting from the vertical axis in the phase-plane analysis of the recent model [A. Constantin, R.S. Johnson, Propagation of very long water waves, with vorticity, over variable depth, with applications to tsunamis, Fluid Dynam. Res. 40 (2008), 175--211]. Our technique can be applied easily in circumstances where the reparametrization device from [A. Constantin, A dynamical systems approach towards isolated vorticity regions for tsunami background states, Arch. Rational Mech. Anal. doi: 10.1007/s00205-010-0347-1] might lead to some serious difficulties.


\noindent{\bf Classification:} 34A12; 34A34; 37E35

\section{Introduction}
A recent model \cite{cj,cja} for studying the currents underlying the water surface before the arrival of large, destructive waves (tsunamis) to the shoreline has been investigated in \cite{carma} via a technical phase-plane analysis relying on several sharp estimates and topological techniques.

A key feature in the long series of arguments there regards the {\it uniqueness\/} of solutions to the following initial value problem
\begin{eqnarray}
\left\{
\begin{array}{ll}
\psi^{\prime\prime}+\frac{\psi^{\prime}}{r}+f(\psi)=0,\quad r\geq r_{0}\geq1,\\
\psi^{\prime}(r_{0})=\psi_{1},\\
\psi(r_{0})=0,
\end{array}
\right.\label{main_ode}
\end{eqnarray}
where $\psi_{1}\neq0$ and the (vorticity) function $f:\mathbb{R}\rightarrow\mathbb{R}$ reads as
\begin{eqnarray}
f(\psi)=\left\{
\begin{array}{ll}
\psi-\frac{\psi}{\sqrt{\vert\psi\vert}},\mbox{ for }\psi\neq0,\\
0,\mbox{ otherwise},
\end{array}
\right.\label{vort1}
\end{eqnarray}
see \cite[Eq. (2.4)]{carma}. Via a reparametrization of the integral curves of the differential equation from (\ref{main_ode}) recast in polar coordinates, the author has established the continuity of the solution with respect to the initial data. We recall that the uniqueness of solution to an initial value problem is essential for establishing the continuity with respect to the initial data when the functional terms of an ordinary differential equation are non-Lipschitzian, cf. \cite[p. 24, Theorem 3.4]{hale}. Constantin's technique, however, might lead to some serious difficulties when one is employing a more complicated (though academic) vorticity function, e.g,
\begin{eqnarray}
f(\psi)=\left\{
\begin{array}{ll}
\psi-\frac{\psi}{\sqrt{\vert \psi\vert}}\left[1+c_{1}-\sin\left(\frac{c_{2}\psi^{2}}{\psi^{2}+1}\right)\right],\quad \psi\neq0,\\
0,\quad \psi=0,
\end{array}
\right.\label{vort2}
\end{eqnarray}
where $0<c_{1}=\sin\frac{c_2}{2}<c_{2}<\frac{3-2\sqrt{2}}{4+3\sqrt{2}}<0.02$, see \cite{must}.

Notice now that both functions from (\ref{vort1}), (\ref{vort2}) satisfy the restrictions
\begin{eqnarray}
\left\{
\begin{array}{ll}
f(0)=0,\quad\psi\cdot f(\psi)<0,\\
\vert f(\psi_{1})-f(\psi_{2})\vert\leq\frac{C}{\sqrt{\min\{\vert\psi_1\vert,\vert\psi_2\vert\}}}\cdot\vert\psi_{1}-\psi_{2}\vert,
\end{array}
\right.\label{main_restr}
\end{eqnarray}
for any numbers $\psi$, $\psi_{1}$, $\psi_{2}$ lying in $[-\delta,0)\cup(0,\delta]$, with $\psi_{1}\cdot\psi_{2}>0$, and some $C$, $\delta>0$.

In the present note, we apply a technique inspired by the approaches from \cite{kk,c97} for dealing with similar issues to prove that the problem (\ref{main_ode}) has a unique solution to the right of $r_0$. Given the preceding formulas of the vorticity functions $f$, it is enough to establish the uniqueness only for the case where $\psi_{1}>0$ (otherwise, we replace $\psi$ with $-\psi$ and remain in the class of nonlinearities $f$ described by (\ref{main_restr})).

\section{Uniqueness of flow}
\begin{theorem}
Assume that the continuous function $f:\mathbb{R}\rightarrow\mathbb{R}$ satisfies the hypotheses (\ref{main_restr}). Then, given $\psi_{1}>0$, the initial value problem (\ref{main_ode}) has a unique solution to the right of $r_0$.
\end{theorem}

\textbf{Proof.} Assume that $\psi$, $\Psi$ are two solutions of (\ref{main_ode}) defined in an interval $[r_{0},r_{1}]$ small enough to ensure that $\psi(r)$, $\Psi(r)\in(0,\delta]$ for every $r\in(r_{0},r_{1}]$. 

Introduce the function $x\in C^{1}([r_{0},r_{1}],\mathbb{R})$ with the formula $x(r)=\Psi(r)-\psi(r)$. Notice as well that $x(r_0)=x^{\prime}(r_0)=0$ and we have the essential estimate of the behavior of $x$ close to $r_0$:
\begin{eqnarray}
\lim\limits_{r\searrow r_{0}}\frac{x(r)}{\ln(r/r_{0})}=\lim\limits_{r\searrow r_{0}}[rx^{\prime}(r)]=r_{0}\cdot x^{\prime}(r_0)=0.\label{estima_nagumo}
\end{eqnarray}
The function $y:[r_{0},r_{1}]\rightarrow[0,+\infty)$ with the formula
\begin{eqnarray*}
y(r)=\left\{
\begin{array}{ll}
\frac{\vert x(r)\vert}{\ln(r/r_{0})},\mbox{ for }r\in(r_{0},r_{1}],\\
0,\mbox{ otherwise},
\end{array}
\right.
\end{eqnarray*}
is continuous by means of (\ref{estima_nagumo}), so, there exists $r^{\star}=r^{\star}(r_{1})\in(r_{0},r_{1}]$ such that $y(r^{\star})=\sup\limits_{s\in(r_{0},r_{1}]}y(s)>0$.

The equation from (\ref{main_ode}) can be recat as
\begin{eqnarray*}
[r\psi^{\prime}(r)]^{\prime}=-rf(\psi),
\end{eqnarray*}
leading to
\begin{eqnarray}
\psi(r)=r_{0}\psi_{1}\ln\frac{r}{r_0}-\int_{r_0}^{r}\tau\ln\frac{r}{\tau}\cdot f(\psi(\tau))d\tau,\quad r\in[r_{0},r_{1}].\label{main_integral}
\end{eqnarray}
The latter integral equation yields
\begin{eqnarray*}
\vert x(r)\vert\leq\int_{r_0}^{r}\tau\ln\frac{r}{\tau}\cdot\frac{C\vert x(\tau)\vert}{\sqrt{\min\{\psi(\tau),\Psi(\tau)\}}}\thinspace d\tau,\quad r\in[r_{0},r_{1}].
\end{eqnarray*}

Given the sign of $f$ in $(0,\delta]$, we get from (\ref{main_integral}) that $\min\{\psi(r),\Psi(r)\}\geq r_{0}\psi_{1}\ln\frac{r}{r_0}$. As a consequence, we obtain that
\begin{eqnarray}
\vert x(r)\vert&\leq&\frac{C}{\sqrt{r_{0}\psi_{1}}}\cdot\int_{r_0}^{r}\tau\ln\frac{r}{\tau}\frac{\vert x(\tau)\vert}{\sqrt{\ln(\tau/r_{0})}}\thinspace d\tau\nonumber\\
&\leq&\frac{C}{\sqrt{r_{0}\psi_{1}}}\cdot\int_{r_0}^{r}\tau\frac{\vert x(\tau)\vert}{\ln(\tau/r_{0})}\thinspace d\tau\cdot\ln\frac{r}{r_0},\quad r\in[r_{0},r_{2}],\label{intermed_1}
\end{eqnarray}
where $r_{2}\in(r_{0},r_{1}]$ is close enough to $r_0$ to have $\ln(r/r_{0})\in[0,1)$ and respectively $\frac{C}{\sqrt{r_{0}\psi_{1}}}\cdot\frac{r^{2}-r_{0}^{2}}{2}\leq\frac{1}{2}$ for all $r\in[r_{0},r_{2}]$.

Further, the inequality (\ref{intermed_1}) reads as
\begin{eqnarray*}
y(r)\leq\frac{C}{\sqrt{r_{0}\psi_{1}}}\int_{r_0}^{r}\tau y(\tau)d\tau<\frac{C}{\sqrt{r_{0}\psi_{1}}}\int_{r_0}^{r}\tau d\tau\cdot y(r^{\star}),\quad r^{\star}=r^{\star}(r_2),
\end{eqnarray*}
for every $r\in(r_{0},r_{2}]$. Finally, we get that
\begin{eqnarray*}
y(r)\leq\frac{1}{2}\cdot y(r^{\star}),\quad r\in(r_{0},r_{2}],
\end{eqnarray*}
which leads to an obvious contradiction by taking $r=r^{\star}$.

The proof is complete. $\square$

\section{Comments}
The reduction of the comparison of two (hypothetical) solutions to an initial value problem to an analysis of the behavior of the solutions to some integral inequality has been employed in various uniqueness results that have stemmed from the classical first order, and respectively $n$-th order Nagumo criteria, see \cite{nag,nag2}. Among the far-reaching generalizations of such results we would like to cite the theorems of Athanassov \cite{ath} devoted to first order equations, and of Wintner \cite{wintner} and Constantin \cite{c97} for the $n$-th order differential equations. A recent reevaluation of the technique has been presented in \cite{cjapan}.

\section{Acknowledgment.}
The author is supported by the FWF--grant I544-N13 {\lq\lq}Lagrangian kinematics of water waves{\rq\rq} of the Austrian Science Fund.

%\label{lastpage}

\end{document}